\documentclass[letterpaper, 10 pt, conference]{ieeeconf}  % Comment this line out if you need a4paper

\IEEEoverridecommandlockouts                              % This command is only needed if 
\usepackage{graphics} % for pdf, bitmapped graphics files
\usepackage{amsmath} % assumes amsmath package installed
\usepackage{amssymb}  % assumes amsmath package installed
\usepackage{cite}
\usepackage[usenames,dvipsnames]{color}%for colours
\usepackage{soul}%to strikethrough

\newtheorem{thm}{Theorem}[section]

\newtheorem{defn}[thm]{Definition}
\newtheorem{eg}[thm]{Example}

\usepackage{graphicx,color,amsmath}

\title{\LARGE \bf
Linearized stability analysis of nonlinear partial differential equations
}

\author{Rasha Al Jamal, Amenda Chow and Kirsten Morris$^\dagger$% <-this % stops a space
\thanks{$^\dagger$ Dept. of Applied Mathematics, University of Waterloo, Waterloo, Ontario N2L 3G1, CANADA,         {\tt\small raljamal@uwaterloo.ca},
   {\tt\small a29chow@uwaterloo.ca},
        {\tt\small kmorris@uwaterloo.ca}}
}

\begin{document}

\maketitle
\thispagestyle{empty}
\pagestyle{empty}

%%%%%%%%%%%%%%%%%%%%%%%%%%%%%%%%%%%%%%%%%%%%%%%%%%%%%%%%%%%%%%%%%%%%%%%%%%%%%%%%
\begin{abstract}
Lyapunov's indirect method is an attractive method for analyzing stability of non-linear systems since only the stability of the corresponding linearized system needs to be determined. Unfortunately, the proof for finite-dimensional systems does not generalize to infinite-dimensions. In this paper a unified approach to  Lyapunov's indirect method for infinite-dimensional system is described. It is shown how existing sufficient conditions fit this framework and a new sufficient condition is presented.  
\end{abstract}

%%%%%%%%%%%%%%%%%%%%%%%%%%%%%%%%%%%%%%%%%%%%%%%%%%%%%%%%%%%%%%%%%%%%%%%%%%%%%%%%

% Motivation
%---------------
\section{INTRODUCTION}

Stability theory for finite-dimensional nonlinear systems is well established \cite{Haddad-2008,Khalil-2002, Sontag-1998}. It is natural to generalize this theory for infinite-dimensional systems. Lyapunov's direct method generalizes to  infinite-dimensional systems in a straightforward manner.  LaSalle's invariance principle also generalizes, provided that the orbit of the system is pre-compact.  \cite{Baker-1969, Luo-1999, Walker-1976, Xu-2003}. However, finding a Lyapunov function for nonlinear infinite-dimensional systems is challenging as is showing the pre-compactness of the orbit of the system. Lyapunov's indirect method is appealing as it is a systematic approach and the theory for the stability of linear systems is well understood. %\achg{\st{Although many researchers use this method to analyze the stability of nonlinear infinite-dimensional systems, they do not provide the justification for this method to work. In particular, there are two obstacles that arise in using Lyapunov's indirect method. First, it is important to determine how to linearize the nonlinear system. That is, a differentiating definition is needed to differentiate the nonlinear operator which is often unbounded. Second, a}} 
However, the proof of Lyapunov's indirect method  for finite-dimensional systems does not generalize to  infinite-dimension dimensions.
% In particular, if the linearized system has some type of stability, how does this determine the stability of the nonlinear system? %\achg{\st{Hans Zwart  determined an example where Lyapunov's indirect method does not hold for an infinite-dimensional system. This example is presented later.}}

Some results on linearization  as a way of analyzing the stability of nonlinear infinite-dimensional systems have been obtained \cite{Henry-1981, Kato-1995, Smoller-1994}. For instance, Smoller \cite[Theorem~11.17]{Smoller-1994} showed that if the nonlinear operator in the model is locally Lipschitz continuous, continuously Fr\'{e}chet differentiable and also a condition related to the nonlinear operator being twice continuously differentiable is satisfied, then the $C_0$-semigroup generated by the nonlinear system is continuously Fr\'{e}chet differentiable. Using this result, Smoller \cite[Theorem~11.22]{Smoller-1994} also showed that if the linearized system at an equilibrium generates an exponentially stable $C_0$-semigroup, then the nonlinear system generates a locally exponentially stable $C_0$-semigroup in a neighbourhood of that equilibrium. Kato \cite[Corollary 2.2]{Kato-1995} relaxed the condition of the nonlinear operator and was able to achieve the same stability result. 
%This result is a case of a more general result Kato \cite[Theorem 20]{Kato-1995} for which the conditions are difficult to verify. 
In \cite{Temam-1988}, the Fr\'{e}chet differentiability of the $C_0$-semigroup corresponding to a class of quasilinear systems is proved under different conditions from the ones in \cite{Kato-1995},  \cite{Smoller-1994}. 

In this paper, %\achg{\st{a more general result on the linearized stability is proved. In particular,}}
it is shown that the Fr\'{e}chet differentiability of the $C_0$-semigroup generated by the nonlinear infinite-dimensional system plays an important role in the justification of Lyapunov's indirect method. Exponential stability as opposed to asymptotic stability is also important.  The nonlinear Kuramoto-Sivashinsky equation is presented as an example. %Lyapunov's indirect method is used to analyze the stability of the nonlinear Kuramoto-Sivashinsky equation with periodic boundary conditions.
A new sufficient condition for  Fr\'{e}chet differentiability is also presented. %\achg{\st{This work is applied to the nonlinear Kuramoto-Sivashinsky equation.}}

%Some preliminaries are introduced, followed in Section 3 by a discussion of the conditions and results needed for  Lyapunov's indirect method  to be used to analyze the local stability of a nonlinear infinite dimensional system. %\achg{\st{a counter example where Lyapunov's indirect method fails to give a stability analysis for a nonlinear system. Then a class of nonlinear systems defined on a Banach space is presented with some conditions under which Lyapunov's indirect method can be used to analyze the stability of an equilibrium.}}
%In particular, it is shown that the Fr\'{e}chet differentiability of the nonlinear $C_0$-semigroup generated by the nonlinear system plays an important role for Lyapunov's indirect method to be applied to infinite-dimensional systems. % Preliminaries
%---------------------
\section{PRELIMINARIES}
\label{prelimanries}

Consider the nonlinear {time-invariant} system defined on a Banach space $X$ with norm $\| \cdot \|$
\begin{eqnarray}
\begin{array}{ll}
\dot{z}\left(t\right)=F\left(z\left(t\right)\right),\ \ \ t\geq 0\\
z\left(0\right) = z_0,
\end{array}
\label{infinite-nonlinear}
\end{eqnarray}
where $z_0$ is the initial condition, the nonlinear operator $F:\mathcal{D}\left(F\right)  \subset X \rightarrow X$ is densely defined on $X$. Assume that this system is well-posed; that is, it has a unique  solution that can be written %\achg{as \sout{in terms of a nonlinear $C_0$-semigroup in $X$}}
\begin{eqnarray*}
z\left(t\right)= S\left(t\right)z_0,
\end{eqnarray*}
where $S\left(t\right)$ is a nonlinear $C_0$-semigroup on $X$ generated by the operator $F$.%\begin{defn}$\left(\right.$Types of stability$\left.\right)$\\ \achg{OMIT?}
% Consider the nonlinear system (\ref{infinite-nonlinear}). 
%\begin{enumerate}
%\item The equilibrium $z_e$ to (\ref{infinite-nonlinear}) is locally stable if for all $\varepsilon>0$, there exists $\delta>0$ such that if $\|z_0-z_e\| <\delta$, then $\|z\left(t\right)-z_e\|<\varepsilon$, $t\geq0$.
%\item The equilibrium $z_e$ to (\ref{infinite-nonlinear}) is locally asymptotically stable if it is stable and there exists $\delta>0$ such that if $\|z_0-z_e\| <\delta$, then $\lim_{t \rightarrow \infty} z\left(t\right)=z_e$.
%\item The equilibrium $z_e$ to (\ref{infinite-nonlinear}) is locally exponentially stable if there exists $\delta, \alpha, \beta >0$ such that if $\|z_0-z_e\| <\delta$, then $\|z\left(t\right)-z_e\| \leq \alpha \|z_0-z_e\| e^{-\beta t}$, $t\geq 0$.
%\item The equilibrium $z_e$ to (\ref{infinite-nonlinear}) is globally asymptotically stable if it is stable and for all $z_0 \in H$, we have $\lim_{t \rightarrow \infty} z\left(t\right)=z_e$.
%\item The equilibrium $z_e$ to (\ref{infinite-nonlinear}) is globally exponentially stable if there exists $\alpha, \beta >0$ such that for all $z_0 \in H$, we have $\|z\left(t\right)-z_e\| \leq \alpha \|z_0-z_e\| e^{-\beta t}$, $t\geq 0$.
%\item The equilibrium $z_e$ to (\ref{infinite-nonlinear}) is unstable if it is not stable.
%\end{enumerate}
%\label{stability-defn}
%\end{defn}

Let $z_e$ is be an equilibrium of (\ref{infinite-nonlinear}); that is,  with $F\left(z_e\right)=0$.  
The standard definitions of stability for an equilibrium point are used here; see for example, \cite{Walker-1976}. Many infinite-dimensional dynamical systems possess an infinite number of stable equilibria. A simple example is the heat equation with  Neumann boundary conditions.

A set of equilibrium points can also be characterized as stable. %\achg{\sout{The equilibrium solution to the dynamical system is determined by the initial condition. } -- Not really needed and also equilibrium solutions also greatly depend on the boundary conditions.}
\begin{defn} \cite[Definition 2.6]{Xu-2003} $\left(\right.$Stable Equilibrium Set$\left.\right)$ \\
%\achg{\st{Consider the dynamical system}(\ref{infinite-nonlinear}).} 
Let $E$ be the set of all equilibria to~(\ref{infinite-nonlinear}). The set E is said to be stable if for every $\varepsilon>0$, there exists $\delta>0$ such that if $dist_{X} \left(z_0 ,E\right)  <\delta$, then
\begin{eqnarray*}
dist_{X}\left(z\left(t\right),E\right) <\varepsilon,\ \ \ t\geq0 
\end{eqnarray*}
\label{stable_set}
where $dist_{X}\left(z,E\right)=inf\{\| z-y \| : y \in E\}$.
\end{defn}

\begin{defn} \cite[Definition 2.6]{Xu-2003} $\left(\right.$Globally Asymptotically Stable Equilibrium Set$\left.\right)$\\
Let $E$ be the set of all equilibria to~(\ref{infinite-nonlinear}). The set E is said to be globally asymptotically stable if it is stable and for every $z_0 \in X$, 
\begin{eqnarray*}
\lim_{t\rightarrow \infty} dist_{X}\left(z\left(t\right),E\right)=0.
\end{eqnarray*}
\label{globally_asymptotically_stable_set}
\end{defn}

% Classes when Lyapunov indirect method holds
%-----------------------------------------------------------------
\section{LINEARIZED STABILITY OF NONLINEAR SYSTEMS}
\label{stability}

%\achg{\st{Analyzing the stability of nonlinear dynamical systems has attracted researchers for many years. One way of analyzing the stability of nonlinear finite-dimensional systems is by Lyapunov's indirect method. That is, the stability analysis of the linearized systems helps in analyzing the stability of the nonlinear system.}} 
%Consider the finite-dimensional dynamical system defined in $\mathbb{R}^n$, for $n<\infty$
%\begin{eqnarray}
%\begin{array}{ll}
%\dot{z}\left(t\right)=f\left(z\left(t\right)\right),\ \ \ t\geq 0\\
%z\left(0\right) = z_0,
%\end{array}
%\label{finite-nonlinear}
%\end{eqnarray}
%where $z_0$ is the initial condition, the function $f:\mathcal{D} \subset \mathbb{R}^n \rightarrow \mathbb{R}^n$ is continuously differentiable and $\mathcal{D}$ is an open set with $0\in \mathcal{D}$. Assume $f\left(z_e\right)= 0$, then $z_e$ is an equilibrium of the nonlinear system. %\achg{\st{Below is the theorem statement for Lyapunov's indirect method for finite-dimensional systems.}}

The following result is well-known.
\begin{thm}(e.g. \cite[Theorem 3.19]{Haddad-2008})\\ 
Consider the nonlinear system  (\ref{infinite-nonlinear}) and assume that $X$ is finite-dimensional. Assume also that $F$ is  differentiable and define
\begin{eqnarray*}
A=\left. \frac{\partial F}{\partial z} \right|_{z=z_e}
\end{eqnarray*}
to be the linearization of (\ref{infinite-nonlinear}). 
Then% \achg{\st{the following statements hold}}:
\begin{enumerate}
\item[(1)] if $Re \lambda<0$ for all $\lambda \in \sigma\left(A\right)$, then the equilibrium $z_e$ to (\ref{infinite-nonlinear}) is exponentially stable where $\sigma\left(A\right)$ is the spectrum of $A$.
\item[(2)] if there exists $\lambda \in \sigma\left(A\right)$ such that $Re \lambda >0$, then the equilibrium $z_e$ to (\ref{infinite-nonlinear}) is unstable.
\end{enumerate}
\label{Lyapunov-indirect-method}
\end{thm}
{The proof of Theorem~\ref{Lyapunov-indirect-method}, see for example,   \cite[Thm~3.1, Thm. 3.12]{Haddad-2008},   relies on showing that a Lyapunov function for the linear system, which can be easily constructed, is also a Lyapunov function  for the non-linear system in a region around the equilibrium point. It follows then that  that $z_e$ is locally exponentially stable.  The proof of instability is similar. 
}

Using Lyapunov's indirect method for nonlinear infinite-dimensional systems requires a justification similar to Theorem \ref{Lyapunov-indirect-method} that the stability of the linearized systems reflects the stability of  the nonlinear system.  However, generalization of the proof to infinite-dimensions is not straightforward since the operator $F$ is typically not Fr\'{e}chet differentiable; in fact it is generally an unbounded operator. There are two issues that need to be addressed. First, how to linearize the nonlinear system defined on a Banach space $X$?  Second, what conditions guarantee that the stability of the linearized infinite-dimensional system is the same as the nonlinear system? That is, if the linearized system is stable or unstable, then does the same conclusion apply to the original nonlinear system?
 
%Consider (\ref{infinite-nonlinear}) as an infinite-dimensional system. %\achg{\st{where $z_e$ is an equilibrium \achg{\st{solution}} to the system. It is of interest to linearize (\ref{infinite-nonlinear}) at the equilibrium $z_e$. }}

\begin{defn} %$\left(\right.$Fr\'{e}chet Differentiable$\left.\right)$
 \cite[Definition 3.1.1]{Lebedev-2002}
Consider an operator $F: X \rightarrow X$ defined on a normed linear space $X$. The operator $F$ is { Fr\'{e}chet differentiable} at $z_0$ if there exists a bounded linear operator $DF\left(z_0\right):X\rightarrow X$ such that  for all $h\in X$
\begin{eqnarray}
\lim_{h\rightarrow 0} \frac{\|F\left(z_0+h\right) - F\left(z_0\right) -DF\left(z_0\right)h  \|}{\|h\|}=0,
\label{frechet-limit}
\end{eqnarray}
That is,
\begin{eqnarray*}
F\left(z_0+h\right) - F\left(z_0\right) = DF\left(z_0\right)h + \omega\left(z_0,h\right),
\end{eqnarray*}
where
\begin{eqnarray*}
\lim_{\|h\| \rightarrow 0}  \frac{\| \omega \left(z_0,h\right)\|}{\|h\|} \rightarrow 0 .
\end{eqnarray*}
The operator $F$ is  said to be Fr\'{e}chet differentiable if it is Fr\'{e}chet differentiable at every $z_0\in X.$
\label{Frechet}
\end{defn}

The next theorem demonstrates that Fr\'{e}chet differentiability of the $C_0$-semigroup generated by the nonlinear infinite-dimensional system plays a key role in the validity of Lyapunov's indirect method.  A similar result was shown in \cite{Desch-1986} under the condition that the number of unstable eigenvalues corresponding to the linearized system is finite. In the next theorem this assumption is not required. 
For more details,  see \cite{Rasha-2013, Rasha-2013c}.
% \achg{Note that the linearized system is the G\^ateaux derivative of the nonlinear system.}
\begin{thm}\label{stability-thm}
%\cite{Rasha-2013,Rasha-2013c}
Consider the nonlinear system (\ref{infinite-nonlinear}) defined on a Banach space $X$. Assume that the nonlinear operator $F:\mathcal{D}\left(F\right) \subset X \rightarrow X$ generates a nonlinear $C_0$-semigroup $S\left(t\right)$. Let $z_e$ be an equilibrium for the above system (\ref{infinite-nonlinear}) and suppose that $S\left(t\right)$ is  Fr\'{e}chet differentiable at $z_e$. 
\begin{itemize}
\item[(i)] If $z_e$ is an exponentially stable equilibrium of the linearized system, then $z_e$ is a locally exponentially stable equilibrium of the nonlinear system (\ref{infinite-nonlinear}).
\item[(ii)] 
If the linearized system is unstable, then the nonlinear system (\ref{infinite-nonlinear}) is locally unstable.
\end{itemize}
\end{thm}
\noindent {\bf Proof:} %Since $z_e$ is an exponentially stable equilibrium solution of the linearized system, then there exists $M\geq 1$ and $\gamma>0$ such that for all $z_0 \in X$
%\begin{eqnarray}
%\|T_{z_e} \left(t\right) z_0 -z_e \| \leq M e^{- \gamma t} \|z_0-z_e\|, \ t \geq 0.
%\label{linear-bound}
%\end{eqnarray}
(i) Let  $ T_{z_e}(t)$ be the Fr\'{e}chet derivative of $S(t)$ at $z_e$. It follows that% (Definition \ref{def-Frechet}) 
\begin{eqnarray*}
S\left(t\right) z_0 - S\left(t\right) z_e = T_{z_e}\left(t\right) \left(z_0-z_e\right) + \omega \left(z_e,z_0-z_e\right),
\end{eqnarray*}
where
\begin{eqnarray}\label{eqlimitforomegaiszero}
\lim_{\|z_0-z_e\| \rightarrow 0}\frac{\|\omega \left(z_e,z_0-z_e \right)\|}{\|z_0- z_e\|}=0.
\end{eqnarray}
That is, for any $t>0$,  $\varepsilon_t >0$, there exists $\delta>0$ such that if $\|z_0-z_e\| <\delta$, 
\begin{eqnarray*}
\frac{\|\omega \left(z_e,z_0-z_e \right)\|}{\|z_0- z_e\|}<\varepsilon_t.
\end{eqnarray*}

Furthermore, since the $C_0$-semigroups $S\left(t\right)$ and $T_{z_e}\left(t\right)$ are continuous in $t$, then the function $\omega$ is continuous in $t$ and for $M\geq 1$ and $\gamma>0$ such that for all $z_0 \in X$
\begin{eqnarray}
\|T_{z_e} \left(t\right) z_0 -z_e \| \leq M e^{- \gamma t} \|z_0-z_e\|, \ t \geq 0,
\label{linear-bound}
\end{eqnarray}
then there exists $\varepsilon>0$, $\bar{t}<\infty$,  such that for $\tau \in [0,\bar{t}\,]$ % $\frac{\|\omega \left(z_e,z_0-z_e \right)\|}{\|z_0- z_e\|}<\varepsilon$ and 
\begin{align}
\|S\left(\tau \right) z_0 - z_e\| &\leq  \| T_{z_e}\left(\tau \right) \left(z_0-z_e\right) \| + \| \omega \left(z_e,z_0-z_e\right)\| \nonumber\\
%&&\nonumber \\
&\leq M e^{-\gamma \tau} \| z_0-z_e\| + \varepsilon \|z_0-z_e\|\nonumber\\
%&&\nonumber\\
%&=& \left(M e^{-\gamma \tau }+ \varepsilon\right)  \| z_0-z_e\| ,\nonumber\\
%&&\nonumber \\
&= C  \| z_0-z_e\| \nonumber
\label{C-bound}
\end{align}
where $C=M + \varepsilon$.
%To show  $z_e$ is locally exponentially stable, to the nonlinear system (\ref{infinite-nonlinear}). 
Choose $\bar{t}= \ln \left(4M\right)/\gamma >0$, then using (\ref{linear-bound})
\begin{eqnarray}
\| T_{z_e}\left(\bar{t}\right) z_0 -z_e\| %&\leq& M e^{-\gamma \bar{t}} \|z_0-z_e\|,\nonumber\\
&\leq& \frac{1}{4} \|z_0-z_e\|.
\label{bound1-k}
\end{eqnarray}
It follows that %Furthermore, using the definition of Fr\'{e}chet differentiable,% (Definition \ref{def-Frechet}),
\begin{align*}
\lim_{\|z-z_e\| \rightarrow 0} &\left\| \frac{S\left(\bar{t}\right) z_0 - S\left(\bar{t}\right) z_e - T_{z_e} \left(\bar{t}\right) z_0 + T_{z_e}\left(\bar{t}\right) z_e }{z_0 -z_e} \right\| \\
=& \lim_{\|z-z_e\| \rightarrow 0} \left\| \frac{S\left(\bar{t}\right) z_0 - T_{z_e} \left(\bar{t}\right) z_0  }{z_0 -z_e} \right\| = 0
\end{align*}
and hence, there exists $\delta>0$ such that if $\|z_0-z_e\| <\delta$, then
\begin{eqnarray}
\| S\left(\bar{t} \right) z_0 - T_{z_e} \left(\bar{t} \right) z_0\| \leq \frac{1}{4} \|z_0 -z_e\|.
\label{bound2-k}
\end{eqnarray}

Using (\ref{bound1-k}) and (\ref{bound2-k}),
\begin{eqnarray}
\| S\left(\bar{t}\right) z_0 - z_e\| %=& \|S\left(\bar{t}\right) z_0 - T_{z_e}\left(\bar{t}\right)z_0 + T_{z_e}\left(\bar{t}\right)z_0 -z_e\|,\nonumber\\
%&\leq&  \|S\left(\bar{t}\right) z_0 - T_{z_e}\left(\bar{t}\right)z_0 \| + \|T_{z_e}\left(\bar{t}\right)z_0 -z_e\|,\nonumber\\
\leq \frac{1}{2} \|z_0-z_e\|
= e^{-\ln2} \|z_0-z_e\|.
\label{bound3-k}
\end{eqnarray}

Let $k>0$ be an integer, then using the semigroup property and (\ref{bound3-k}),
\begin{eqnarray}
\|S\left(k\bar{t}\right) z_0 - z_e\|  &=& \| S^k \left(\bar{t}\right) z_0 -z_e\|\nonumber\\
&=& \|S\left(\bar{t}\right) S^{k-1}\left(\bar{t}\right) z_0 -z_e\|\nonumber\\
&\leq& e^{-\ln2} \| S^{k-1}\left(\bar{t}\right) z_0 -z_e\|\nonumber\\
&\leq& e^{-\left(\ln2\right) k} \|z_0-z_e\|.
\label{e-bound}
\end{eqnarray}

For $t>0$, let $k=\left[ t/\bar{t}\,\right]$ and $\tau= t-k\bar{t}$. Then $\tau \in [0,\bar{t}]$ and using the semigroup property, (\ref{bound1-k}) and (\ref{e-bound}),
\begin{eqnarray*}
\|S\left(t\right)z_0 - z_e\| &=& \|S\left(k\bar{t} + \tau\right)z_0 - z_e \|\\
&=& \|S\left(\tau\right) S\left(k\bar{t}\right) z_0 - z_e\|\\
&\leq& C\| S\left(k\bar{t}\right)z_0 -z_e\|\\
&\leq& C e^{-\left(\ln2\right) k} \|z_0-z_e\|\\
&\leq& Ce^{-\alpha t} \|z_0-z_e\|
\end{eqnarray*}
for $\alpha \leq \ln 2/\bar t$. This implies that the equilibrium $z_e$ to the nonlinear system is locally exponentially stable.

(ii) The result is shown by proving the contrapositive. Let $z_e$ be a locally stable equilibrium to the nonlinear system (\ref{infinite-nonlinear}). %Using the definition of Fr\'{e}chet differentiable (Definition (\ref{def-Frechet})), there is $r>0$ so that for all $z_0$ \begin{eqnarray*}
%S\left(t\right)z_0 - S\left(t\right)z_e = T_{z_e}\left(t\right) \left(z_0-z_e\right) + \omega \left(z_e,z_0-z_e\right),
%\end{eqnarray*}
%where $\omega \left(z_e,z_0-z_e\right)$ satisfies
%\begin{eqnarray}
%\lim_{\|z_0-z_e\| \rightarrow 0} \frac{\|\omega \left(z_e,z_0-z_e\right)\|}{\|z_0-z_e\|}=0.
%\label{omega1}
%\end{eqnarray}
Since $T_{z_e}\left(t\right)$ is a linear operator and $z_e$ is an equilibrium,
\begin{eqnarray}
S\left(t\right)z_0 - z_e = T_{z_e}\left(t\right) z_0-z_e + \omega \left(z_e,z_0-z_e\right).
\label{omega3}
\end{eqnarray}

The definition of locally stable equilibrium of the nonlinear system implies that for any $\varepsilon>0$, there exists $\delta>0$ such that if $\|z_0-z_e\| < \delta,$
%\begin{eqnarray}
%\|z_0-z_e\| < \delta,
%\label{delta1}
%\end{eqnarray}
then
\begin{eqnarray*}
\|S\left(t\right)z_0 -z_e\| \leq \frac{\varepsilon}{2}, \textup{    for all } t \geq0,
%\label{delta2}
\end{eqnarray*}

%Also, since
%\begin{eqnarray}
%\lim_{\|z_0-z_e\| \rightarrow 0} \frac{\|\omega \left(z_e,z_0-z_e\right)\|}{\|z_0-z_e\|}=0,
%\label{delta3}
%\end{eqnarray}
From (\ref{eqlimitforomegaiszero}), there is $\hat{\delta}$, with $0<\hat{\delta}<\delta$, such that if $\|z_0-z_e\| \leq \hat{\delta}$, then
\begin{eqnarray*}
\frac{\|\omega \left(z_e,z_0-z_e\right)\| }{\|z_0-z_e\|} \leq \frac{\varepsilon}{2}.
\end{eqnarray*}

Then, from (\ref{omega3})
\begin{align*}
\| T_{z_e}\left(t\right) z_0-z_e \|& \leq \| \omega \left(z_e,z_0-z_e\right)\| +   \| S\left(t\right)z_0 - z_e\|\\ 
&\leq \varepsilon.\,\, 
\end{align*}
Thus, $z_e$ is a stable equilibrium point of the linearization. $\Box$
%\chg{Put in proof from paper with Rasha; edited if necessary so paper is 6 pages. } 
 %For a complete proof see \cite{Rasha-2013c}.

%\begin{thm}
 %Let $z_e$ be an equilibrium for the system (\ref{infinite-nonlinear}). Assume that $S\left(t\right)$ is  Fr\'{e}chet differentiable at $z_e$ and the derivative is given by $T_{z_e}\left(t\right)$. If  \end{thm}
%\achg{It turns out that t}he Fr\'{e}chet differentiability of the semigroup of~(\ref{infinite-nonlinear}) is \achg{crucial in showing the stability of the linearized system is the same as the associated nonlinear system.}  Given this result and assuming exponential stability of the linearized system, it \achg{turns out} that the nonlinear system is locally exponentially stable.
% Example by Hans:
%--------------------------
The requirement in Theorem~\ref{stability-thm} that the linear system exhibits exponential stability is crucial. %\st{ However, it should be noted that Lyapunov's indirect method is not always helpful when used to analyze the stability of nonlinear infinite-dimensional systems.}} 
Below is an example due to  Hans Zwart \cite{Hans} illustrating this point.  The example also highlights a fundamental difference between finite and infinite-dimensions; that is, exponential and asymptotic stability are not equivalent for linear systems in infinite-dimensions.
 
\begin{eg}  \label{example-Hans-Zwart}
%\chg{edit if necessary so paper is 6 pages}
Let $\ell_2$ be the space of square summable sequences and $\mathbb{N}$ the set of natural numbers with norm $||\cdot||_{\ell_2}$. For any $z(t)=(z_1(t),z_2(t),\dots,z_n(t),\dots) \in \ell_2$ with $n\in \mathbb{N}$, consider 
\begin{equation}\label{zwart-system}
\dot{z}_n=-\frac{1}{n} z_n+z_n^2.
\end{equation}
This system has infinitely many equilibrium $z_e\in \ell^2$ since $\dot{z}_n = 0$ if and only if $-\frac{1}{n} z_n + z_n^2 =0$ for $n\in\mathbb{N}$. This implies that $z_n =0,\frac{1}{n}$.  Therefore, the set of equilibria is 
\begin{eqnarray*}
E = \left\{ z\in \ell^2 |  \; z_{n} \in \left\{ 0,\frac{1}{n}\right\}, n\in\mathbb{N} \right\}.
\end{eqnarray*}

%Consider the stability of the zero equilibrium solution $z_e=0$.

Linearize the system (\ref{zwart-system}) around  $z_e=0$ to obtain
\begin{eqnarray}
\begin{array}{ll}
\dot{z}_n \left(t\right) =-\displaystyle{\frac{1}{n}}z_n\left(t\right),\ t\geq0\\
%z\left(0\right) = z_0,
\end{array}
\label{zwart-system-linear}
\end{eqnarray}
which has solution
\[
z(t)=(z_1(0)e^{-t}, z_2(0)e^{-\frac{1}{2}t},\dots).
\]

The linearized system (\ref{zwart-system-linear}) is asymptotically stable since
\begin{eqnarray*}
\lim_{t\rightarrow \infty} \|z\left(t\right) - z_e \|_{\ell_2} &=& \lim_{t\rightarrow \infty} \|z\left(t\right)\|_{\ell_2}\\
&=&  \lim_{t\rightarrow \infty} \left(\sum_{n=1}^{\infty} z_{n}^2\left(0\right) e^{-\frac{2}{n} t} \right)^\frac{1}{2}\\
&=&0.
\end{eqnarray*}

Now consider the stability of the original nonlinear system (\ref{zwart-system}).
The solution to (\ref{zwart-system}) is 
\begin{eqnarray}
z_n \left(t\right) = \frac{z_{0n} e^{-\frac{1}{n} t} }{z_{0n} n (-1 + e^{-\frac{1}{n} t} ) +1 }
\label{exact-sol-zwart}
\end{eqnarray}
where $z_{0n}$ is the initial condition.
For any $\delta>0$, choose $n$ such that $\frac{1}{n} < \delta$. 
In the nonlinear system (\ref{zwart-system}), choose components of the initial condition $z_0$ to be zero except in the $n^{th}$ position, which is chosen to be $\frac{1}{n}$; that is,
\begin{eqnarray*}
z_0 = \left( 0,\cdots,0,\frac{1}{n},0,\cdots\right).
\end{eqnarray*}
Given this initial condition, the solution to (\ref{zwart-system}) is
\begin{eqnarray*}
z\left(t\right) =\begin{array}{c}\displaystyle{\left(0, \cdots, \frac{1}{n}, \cdots \right)} \end{array} 
\end{eqnarray*}
and hence $\|z_0 - z_e \|_{\ell_2} = \frac{1}{n} <\delta$. However,
\begin{eqnarray*}
\lim_{t\rightarrow \infty} \|z\left(t\right) - z_e \|_{\ell_2} = \frac{1}{n} \neq 0.
\end{eqnarray*}
Hence, the zero equilibrium $z_e$ to the nonlinear system (\ref{zwart-system}) is not asymptotically stable.
%\achg{\st{Use of  Lyapunov's indirect method for nonlinear infinite-dimensional systems requires justification that the stability of the linearized systems reflects the stability of  the nonlinear system. The proof for finite-dimensional systems does not generalize to infinite-dimensions. This is clearly demonstrated by Example}~\ref{example-Hans-Zwart}.}

Note that if the solution in~(\ref{exact-sol-zwart}) is truncated to $N$ dimensions 
\begin{align*}
\lim_{t\rightarrow \infty} ||z(t)||^2_{\ell_2} &= \lim_{t\rightarrow \infty}\sum_{n=1}^N \left| \frac{z_{0n} e^{-\frac{1}{n} t} }{z_{0n} n (-1 + e^{-\frac{1}{n} t} ) +1 }\right|^2\\
&=\sum_{n=1}^N\lim_{t\rightarrow \infty} \left|  \frac{z_{0n} e^{-\frac{1}{n} t} }{z_{0n} n (-1 + e^{-\frac{1}{n} t} ) +1 }\right|^2=0,
\end{align*}
and hence the approximated solution is asymptotically stable. The lack of stability is only apparent with the exact solution, not with the approximated solution.$\ \ \Box \\$
\end{eg}

Since it is generally difficult to obtain a  closed form representation of the semigroup,  linearization and stability analysis is generally done using the generator $F$ in (\ref{infinite-nonlinear}). It is desirable to have conditions for  linearized stability in terms of the generator. If the generator is Fr\'{e}chet  differentiable, it is not difficult to show that the semigroup generated by the linearization corresponds to the linearization of the original semigroup and Theorem \ref{stability-thm} can be used.
However,  the generator in an infinite-dimensional space is typically unbounded and in these cases the generator is not   Fr\'{e}chet differentiable.   G\^{a}teaux differentiability is generally a more useful concept for  linearization of the generator $F$.

 \begin{defn}
 %$\left(\right.$G\^{a}teaux Differentiable$\left.\right)$ \\% \achg{\cite[Definition 3.1.2]{Lebedev-2002}}\\
 Let $F:\mathcal{D}\left(F\right) \subset X \rightarrow X$ be an operator defined on a Banach space $X$. The operator $F$ is {G\^{a}teaux differentiable}  at $z_0\in \mathcal{D}\left(F\right)$ if there exists a linear operator $dF\left(z_0\right): X \rightarrow X$ such that
\begin{eqnarray*}
\lim_{\varepsilon \rightarrow 0} \frac{F\left(z_0 +\varepsilon h \right) -F\left(z_0\right)}{\varepsilon}=dF\left(z_0\right) h,
\end{eqnarray*}
where $h, \left(z_0+ \varepsilon h\right) \in \mathcal{D}\left(F\right)$.
\label{def-gateaux}
\end{defn}

The question is then, once the generator is linearized via  a  G\^{a}teaux derivative, does the linearization generate a semigroup; and if so, does this semigroup correspond to the  Fr\'{e}chet derivative of the original system?

The situation for quasilinear systems is fairly well-understood. %\achg{\st{It is useful to investigate when Lyapunov's indirect method can be used to give information about the stability of a nonlinear infinite-dimensional system. There are some results justifying the use of linearization to analyze the stability of  a nonlinear infinite-dimensional system. Often exponential stability of the linearized system is required.}} 
Consider a {time-invariant} quasilinear  system on a Banach space $X,$ %with norm $||\cdot||_Z$ and inner product $\langle \cdot,\cdot\rangle_Z$ 
\begin{eqnarray}
\begin{array}{l}\label{eqstatespaceformGeneral}
\dot{z}\left(t\right) = Az\left(t\right) + f\left(z\left(t\right)\right),\\
z\left(0\right) = z_0,
\end{array}
\end{eqnarray}
where $z\left(t\right) \in X$ is the state and $z_0$ is the initial condition. The operator $A:\mathcal{D}\left(A\right) \subset X \rightarrow X$ is a linear operator that generates a $C_0$-semigroup on $X$ and the nonlinear operator $f:\mathcal{D}\left(f\right) \subset X\rightarrow X$ is Fr\'{e}chet differentiable with  $Df(z)$  the Fr\'echet derivative of $f$ at $z$. It is straightforward to show that $A+ Df(z)$ is the G\^{a}teaux derivative of $Az+f(z)$ at $z$. The linearized system corresponding to  (\ref{eqstatespaceformGeneral})  at the equilibrium point  $z\in Z$ is
\begin{equation}\label{eqstatespaceformGeneralLINEAR}
\frac{d\psi}{dt}=A\psi+Df(z)\psi
\end{equation}

%\begin{thm}
%  \cite[Theorem 11.18]{Smoller-1994}
 %  Smoller assumes $f$ is continuously Fr\'echet differentiable and $f$ is locally Lipschitz continuous, and he also requires that $f$ satisfies
%\[
%||f(z)-f(w)-df_w(z-w)|| \leq c||z-w||^2
%\] for all $z,w$ in a bounded set; given these conditions, the nonlinear $C_0$-semigroup corresponding to ~(\ref{eqstatespaceformGeneral}) is continuously Fr\'{e}chet differentiable
%\chg{quote this theorem formally}
%\end{thm}
  
 The following theorem is a special case of the more general result in \cite[Thm. 2.1]{Kato-1995} for which the conditions are difficult to check. 
 This theorem  generalizes an earlier result  \cite[Theorem 11.22]{Smoller-1994}
 which has more restrictive conditions on $f$. %\achg{Rasha please check reference}

 \begin{thm}Consider equation~(\ref{eqstatespaceformGeneral}). Suppose $A$ generates a $C_0$-semigroup and $f$ is Fr\'echet differentiable on $X$, and that the Fr\'echet derivative of $f$  satisfies
 \[
 ||Df(z_1)-Df(z_2)||\leq c(r)||z_1-z_2||,% \quad \forall ||z_1||\leq r,||z_2||\leq r
 \]
for some $r>0$ and for all $ ||z_1||\leq r,||z_2||\leq r,$ where $c:[0,\infty)\rightarrow [0,\infty)$ is a continuous increasing function. Let $z_e$ be an equilibrium point of (\ref{eqstatespaceformGeneral}).
%and assume  that $A+Df(z_e)$ generates a $C_0$-semigroup on $X . 
If $A+Df(z_e)$ generates an exponentially stable semigroup, then $z_e$ is a locally exponentially stable equilibrium point.
 Conversely, if the linearization is unstable, the original system is also unstable.  
 \end{thm}
 
 \noindent
 {\bf Proof:}
 This is essentially shown in \cite{Kato-1995}. In section 3 of that paper, it is shown that the assumptions imply that the nonlinear semigroup is Fr\'echet differentiable at any equilibrium $z_e$, with generator $A+ dF(z_e) .$ In \cite[Cor. 2.2]{Kato-1995}  it is then shown that exponential stability of the linear semigroup implies local exponential stability of the original system, or Theorem \ref{stability-thm} can be used. Theorem \ref{stability-thm}  implies instability of the original system if the linearization is unstable. 
 
% Thus,  if the system linearized at an equilibrium is exponentially stable, then the nonlinear system is locally exponentially stable \cite{Kato-1995} or Theorem \ref{stability-thm}.  
% \chg{I think Kato showed instability as well?}
 
 %Similar results to Smoller \cite[Theorem 11.22]{Smoller-1994} and Kato \cite[Cor. 2.2]{Kato-1995} are also found
% Furthermore, Temam \cite[Section VI.8]{Temam-1988} considered a more general class of nonlinear systems defined on $X=Z$, where $Z$ is a Hilbert space. but for different class of nonlinear systems. \rchg{Fill} conditions on $A$ and $f$ in (\ref{eqstatespaceformGeneral}). Rasha -- please put details about TEMAM here.

%, Kato relaxed the condition on the nonlinear operator $f$, where he only requires the nonlinear operator $f$ to be continuously Fr\'{e}chet differentiable in a weaker sense, where the uniform convergence of the difference between the operator an its derivative is not required. 

%It is worth mentioning that the above conditions are not the only conditions to achieve a Fr\'{e}chet differentiable $C_0$-semigroup. %\achg{moved sentence to the end of paragraph -- \sout{For instance, in} \cite[Section VI.8]{Temam-1988} \sout{Temam showed the differentiability of the $C_0$-semigroup corresponding to a class of quasilinear systems where the nonlinear operator need not to be bounded.}} \achg{\sout{Moreover,} 

The assumptions on $f$ in  the following theorem are slightly different to those above. 
%similar to \cite[Theorem 11.18]{Smoller-1994}; however, the requirements on $f$ are different. %\achg{Rasha please check reference.}
\begin{thm}
\label{thmlinearizedsemigroup}
Let $Z$ be a Hilbert space with norm $||\cdot||_Z$ and inner product $\langle \cdot,\cdot\rangle_Z$. Consider the quasilinear equation in (\ref{eqstatespaceformGeneral}) and suppose it generates a semigroup, $S(t)$. Assume $\mathrm{Re}\langle Az,z\rangle_{Z} \leq 0$ for all $z \in N_{z,r}=\{ p \in Z: ||p-z||_Z\leq r\}$ %\achg{.\sout{ and that $A$ generates a \rchg{$C_0$-}semigroup.} 
and suppose $f$ is Fr\'echet differentiable on $N_{z,r}$ and its derivative, $Df,$ is locally Lipschitz continuous on $N_{z,r}$. %; that is,
%\[
%||Df_y-Df_z||_Z\leq L_{z,r}||y-z||_Z \mbox{ for all } y \in N_{z,r}
%\]
%where $L_{z,r}$ is the Lipschitz constant that depends on $z$ and some $r>0$. 
Also, for some positive constant $K_{z,r}$ that depends on $z$ and $r$, assume that 
\[
\sup_{\eta \in N_{z,r}}||Df(\eta)||_{op} =K_{z,r}<\infty
\]
where $||\cdot||_{op}$ is the operator norm. Then (\ref{eqstatespaceformGeneralLINEAR}) generates the semigroup, $T_{z}(t),$ %: {Z}\rightarrow {Z}$ 
and %\achg{\sout{furthermore,}}
\[
T_{{z}}(t)=DS(z_0)(t)
\]
where $DS(z_0)(t)$ is the Fr\'{e}chet derivative of $S(t)$ at ${z(0)=z_0}$  for all $0\leq t\leq t_f$ for some positive $t_f$.
\end{thm}

\noindent {\bf Proof:}
Since $A$ generates a semigroup and $Df(z)$ is bounded, then $A+Df(z)$ generates a semigroup, $T_z(t),$ and $\psi(t)=T_z(t)z_0$ is the unique solution to (\ref{eqstatespaceformGeneralLINEAR}) \cite[Theorem~3.2.1]{Curtain-1995}.

Let  $y, z \in N_{z,r}$ be solutions to (\ref{eqstatespaceformGeneral}) and define $h(t):=y(t)-z(t)$.  Taking the time derivative of $h$ and then the inner product with $h$ leads to 
\[
\frac{1}{2}\frac{d}{dt}||h||_Z^2=\mathrm{Re}\langle Ah, h\rangle_Z+\mathrm{Re}\langle f(y)-f(z),h \rangle_Z.
\]
Since $\mathrm{Re}\langle Ah, h\rangle_Z \leq 0$ and applying the Cauchy-Schwarz inequality %(Curtain and Zwart \cite[Definition~A.2.24b]{Curtain1995})
\begin{equation}\label{eqinequalityODE}
\frac{1}{2}\frac{d}{dt}||h||^2_Z \leq || f(y)-f( z)||_Z||h||_Z.
\end{equation}
From the Mean Value Theorem,
\[
||f(y)-f(z)||_Z \leq\sup_{z\in N_{z,r}}|| D f (z) ||_{op}||y-z||_Z.
\]
and hence
\begin{equation*}%\label{eqLipschitzf}
||f(y)-f(z)||_Z \leq K_{z,r} ||y-z||_Z \mbox{ for all } y,z \in N_{z,r}
\end{equation*}
since $K_{z,r}=\sup_{\eta \, \in \, N_{z,r}} ||D f (\eta ) ||_{op}<\infty$. Equation~(\ref{eqinequalityODE}) then becomes
\[
\frac{d}{dt}||h||_Z^2 \leq 2K_{z,r}||h||_Z^2.
\]
Integrating with respect to $t$ yields
\begin{equation}\label{eqexponentialinequality}
||y(t)-z(t)||_Z^2 \leq ||y(0)-z(0)||^2_Ze^{2K_{z,r}t}.
\end{equation}
%Substituting in $h(t)=y(t)-z(t)$ leads to the desired result.

Define $\phi(t):=y(t)- z(t)-\psi(t)$ where $\psi(t)$ is the solution to (\ref{eqstatespaceformGeneralLINEAR}). Taking the time derivative of $\phi$ and then the inner product with $\phi$ leads to
\begin{align*}
\frac{1}{2}\frac{d}{dt}||\phi||^2_Z&=\mathrm{Re}\langle A\phi,\phi \rangle_Z+\mathrm{Re}\langle D f (z) \phi, \phi \rangle_Z\\
&+\mathrm{Re}\langle f(y)-f(z)- D f (z) (y-z), \phi \rangle_Z.
% \leq ||h||_\mathrm{X}||\phi||_\mathrm{X} +||B^\prime(t,x)||_\mathrm{X}||\phi||_\mathrm{X}^2
\end{align*}
Since $f$ is Fr\'echet differentiable, then $||Df(z)||_{op}$ is bounded by a positive constant, $M_z$, and since the derivative of $f$ is locally Lipschitz continuous, and $\mathrm{Re}\langle Az,z\rangle_{Z} \leq 0$, then
\[
\frac{1}{2}\frac{d}{dt}||\phi||^2_Z \leq M_z||\phi||_Z^2+  \frac{L_{z,r}}{2}  ||y - z ||_Z^2||\phi||_Z
\] 
where $L_{z,r} $ is the Lipschitz constant. Applying Gronwall's inequality \cite[Lemma~2.8]{Robinson-2001} with $\phi(0)=0$ implies 
\[
||\phi||_Z \leq L_{z,r}e^{2M_zt}\int_0^te^{-2M_zs}|| y(s)- z(s)||^2_Zds.
\] 
Applying  equation~(\ref{eqexponentialinequality}) for $t\in [0,t_f]$ with $t_f$ any positive constant,
\[
||\phi||_Z \leq L_{z,r}e^{2M_zt_f}|| y_0- z_0||^2_Z\int_0^{t_f}e^{2(K_{z,r}-M_z)s}ds.
\] 
and solving the integrals leads to 
\[
||\phi(t)||_Z \leq k(t_f)|| y_0-z_0||^2_Z, \qquad \mbox{ for } t\in[0,t_f]
\] 
where
\[
k(t_f)=\frac{L_{z,r}}{2\left(K_{z,r}-M_z\right)}(e^{2K_{z,r}t_f}-e^{2Mt_f}).
\]
%We can assume without lose of generality that $k(t_f)$ is well-defined. 
It follows that 
\[
||y- z-\psi||_Z\leq  k(t_f) || y_0- z_0||_Z^2
\]
and hence %From (\ref{eqsemigroupT}) and recalling that $y(t)=F(t)y_0$ and $ z(t) = F(t)z_0$, we obtain% since $\bar z$ is an equilibrium of (\ref{eqstatespaceform}), 
\[
||S(t)y_0-S(t)z_0-T_{z}(t)\psi_0||_Z\leq  k(t_f) || y_0- z_0||_Z^2.
\]
Defining $h_0= y_0 -z_0$ leads to
%\[
%||F(t)(h_0+z_0)-F(t)z_0-T_{ z}(t)h_0||_Z\leq  k(t_f) ||h_0||_Z^2
%\]
%and hence
\[
\lim_{||h_0||_Z\rightarrow 0}\frac{||S(t)(h_0+z_0)-S(t)z_0-T_{ z}(t)h_0||_Z}{||h_0||_Z} = 0
\]
for all $t\in [0,t_f]$.
More details can be found in \cite[Theorem~2.23]{ChowThesis2013}. $\Box$%Therefore, (\ref{eqFrechet}) is satisfied which means $F(t)$ is Fr\'{e}chet differentiable at $ z_0$ and its derivative is $dF_{z_0}(t)=T_{z}(t)$ for $t\in[0,t_f].$
%\chg{put in proof, won't appear anywhere else}

Some limited results have  been achieved for systems that are not quasilinear.
\begin{thm}
\cite[Section VI.8]{Temam-1988} 
\label{temam}
%\chg{State result formally}
Let $Z$ be a Hilbert space with norm $||\cdot||_Z$ and inner product $\langle \cdot,\cdot\rangle_Z$. Consider equation~(\ref{eqstatespaceformGeneral}) with $A:D(A)\rightarrow Z$ closed, negative and self-adjoint. Define $Y=D((-A)^{1/2})$ with norm   $||y||_Y=||(-A)^{1/2}y||_Z$  and  the dual space $Y^\prime$  with norm $||y||_{Y^\prime}=||(-A)^{-1/2}y||_Z$.
Assume that
\begin{equation*}%\label{eqTemamCondition1}
 f(z(t))-f(w(t)) = L(z(t)-w(t))+Q(z(t)-w(t))
\end{equation*}
where $L$ is a linear bounded operator on $Y$ to $Y^\prime$ such that for  some  $0<\epsilon \leq 1$ and positive constant $c_\epsilon$ that depends on $\epsilon$,
\begin{equation*}%\label{eqTemamCondition2}
|\langle Lv,v\rangle_Z|\leq (1-\epsilon)||y||^2_\mathrm{Y}+c_\epsilon||y||_Z^2
\end{equation*}
for all $y\in Y$ and also assume $Q$ satisfies
\begin{equation*}%\label{eqTemamCondition3}
||Q(z(t)-w(t))||_{\mathrm{Y}^\prime}\leq k_1||z-w||_\mathrm{Y}^{1+\sigma_1}
\end{equation*}
for some $k_1>0$ and $\sigma_1>0$.   Also assume that for every $R>0$  there exists $0< \sigma_2\leq 1$ and constant ${k}_R$ depending on $R$ such that
\[
|\langle f(z) -f(w) ,z-w\rangle_Z|\leq k_R ||z-w||_Z^{\sigma_2}||z-w||_{V}^{2-\sigma_2}
\]
for all $z,w \in Y$ with $||z||_Z\leq R$ and $||w||_Z\leq R$.  
%and the nonlinear operator $f$ satisfies
%\begin{eqnarray*}
%\left| \langle f\left(z\right) - f\left(w\right) , z-w \rangle_Z \right| \leq c \|z-w\|^{\eta_0}_{Z} \|z-w\|^{2-\eta_0}_V,
%\end{eqnarray*}
%for $z,w \in V \subset Z$, $0<\eta_0 \leq 1$ and $z,w$ being in a bounded set. Moreover, the operator $f$ can be written as $  f\left(z\right) - f\left(w\right) =l_0 \left(z-w\right) + l_1\left(z-w\right)$ where $l_0$ is a linear bounded operator that satisfies $\langle l_0 z , z \rangle_Z \leq \left(1-\varepsilon\right) \|z\|_V^2 + C_\varepsilon \|z\|_Z^2$, where $0<\varepsilon \leq 1$ and there exists $\eta_1>0$ and $c>0$ such that $\| l_1\left(z-w\right)\|_{-1} \leq c \|z-w\|_V^{1+\eta_1}$. 
Given these conditions, the semigroup of (\ref{eqstatespaceformGeneral}) is Fr\'echet differentiable at any $z$  with its derivative equal to the semigroup generated by $A+ d f (z)  . $
\end{thm}

Examples that satisfy the assumptions of Theorem~\ref{temam} are found in \cite{Temam-1988}. These include special cases of the Navier-Stokes  and wave equations.

The above classes are not exhaustive. Consider for example the Kuramoto-Sivashinsky (KS) equation \cite{Kuramoto-1978,Sivashinsky-1977,Temam-1988} 
%The Kuramoto-Sivashinsky (KS) equation is a nonlinear partial differential equation that describes reaction-diffusion models and chaotic phenomena such as those that occur in turbulence\cite{Armaou-2000, Dubljevic-2010, Gustafsson-2010, Koboyashi-2002, Lee-2005, Liu-2001, Sakthivel-2006}.  %\achg{\st{that is first-order in time and fourth-order in space. It models reaction-diffusion systems and is related to various pattern formation phenomena where turbulence or chaos appear. For instance, it models long wave motions of the liquid film over a vertical plane, dendritic fronts in dilute binary alloys, unstable flame front, Belouzov-Zabotinskii reaction pattern and interfacial instabilities between two viscous fluids} \cite{Armaou-2000, Dubljevic-2010, Gustafsson-2010, Koboyashi-2002, Lee-2005, Liu-2001, Sakthivel-2006}.} \achg{\st{Also, it describes the nonlinear saturation mechanism of dissipative trapped ion modes} \cite{Liu-2001}. }
%Furthermore, the KS equation model is often used in the study of convective hydrodynamics, plasma confinement in toroidal devices, interfacial instabilities between two viscous fluids \cite{Dubljevic-2010}, the bifurcation solutions of the Navier-Stokes equation \cite{Sell-2000} and the B\'{e}nard problem in an elongated box. The convection cells patterns developed from heating the plane horizon from below are modelled by the nonlinear KS equation \cite{Temam-1988}.
with periodic boundary conditions defined on the Hilbert space $L^2(-\pi,\pi)$
\begin{eqnarray}
\begin{array}{l}
\displaystyle{\frac{\partial z}{\partial t} + \nu \frac{\partial^4 z}{\partial x^4} + \frac{\partial^2 z}{\partial x^2} + z \frac{\partial z}{\partial x} =0},\ \ t\geq 0,\\
\ \\
\displaystyle{\frac{\partial^n z}{\partial x^n}\left(-\pi,t\right) = \frac{\partial^n z}{\partial x^n}\left(\pi,t\right)},\ \ n=0,1,2,3,\\
\ \\
z\left(x,0\right)= z_0\left(x\right),
\end{array}
\label{ks}
\end{eqnarray}
where $z\in L^2(-\pi,\pi)$ is the state of the system, $\nu>0$ is the instability parameter, $-\nu \frac{\partial^4 z}{\partial x^4}$ is the dissipative term, $\frac{\partial^2 z}{\partial x^2}$ is the anti-dissipative term and $ z \frac{\partial z}{\partial x}$ is the nonlinear term \cite{Sell-2000}. %\achg{\st{that can be interpreted as the energy transfer mechanism that transfers energy from low to high wave numbers and hence prevents the low wave numbers from blowing up} \cite{Sell-2000}. } 
This equation has a unique strong solution $z\left(t\right) = S\left(t\right) z_0$, where
\begin{eqnarray*}
z\left(t\right) \in L^2\left([0,T]; H^2_{per}(-\pi,\pi) \right) \cap L^\infty \left([0,T];L^2(-\pi,\pi)\right),  
\end{eqnarray*}
and $S\left(t\right)$ is a nonlinear $C_0$-semigroup  \cite[Theorem 5.4.3]{Sell-2000}. The stability analysis of the KS equation depends on the parameter $\nu$. If the instability parameter $\nu>1$,  the set of all constant equilibria is globally asymptotically stable. Furthermore, if $\nu=1$, then the zero equilibrium is Lyapunov stable. This is proven using a  Lyapunov function and LaSalle's invariance principle \cite{Rasha-2013,Rasha-2013c},

%Many researchers studied the stability of the dynamics of the KS equation analytically as well as numerically \cite{Armbruster-1989, Brown-1992, Byrnes-2010, Chen-1986,  Christofides-1998c, Coward-1995, Elgin-1996, Foias-1989, Guo-1996, Hyman-1986, Kevrekidis-1990, Lou-2003, Nicolaenko-1985}. The KS equation is a nonlinear PDE that can be regarded as an extension of the heat equation \cite{Liu-2001}. Furthermore, analytical as well as numerical studies on the dynamics of the KS equation showed the existence of unstable steady-state and periodic solutions and chaotic behaviour for very small values of the instability parameter \cite{Lou-2003}. In \cite{Zhang-2011}, Zhang, Song and Axia considered the KS equation with periodic boundary conditions and odd solutions and were able to show that the zero equilibrium solution to the KS equation is globally exponentially stable for certain values of the instability parameter. A more general result is obtained in this section.

Define  $A: \mathcal{D}\left(A\right) = H^4_{per}(-\pi,\pi) \subset H^4 (-\pi,\pi) \rightarrow L^2(-\pi,\pi)$ by
\begin{eqnarray}
Az=-\nu \frac{\partial^4 z}{\partial x^4} - \frac{\partial^2 z}{\partial x^2}
\label{A}
\end{eqnarray}
and the nonlinear operator $J:\mathcal{D}\left(J\right)=H^1_{per} (-\pi,\pi)\subset H^1 (-\pi,\pi) \rightarrow L^2(-\pi,\pi)$ by
\begin{eqnarray}
J\left(z\right) = -z \frac{\partial z}{\partial x}.
\label{J}
\end{eqnarray}
The KS equation (\ref{ks}) can be written 
\begin{eqnarray}
\begin{array}{l}
\dot{z} = Az + J\left(z\right), \\
z\left(0\right)= z_0.
\end{array}
\label{cauchy-ks}
\end{eqnarray}
The G\^{a}teaux derivative $dJ: H^1(-\pi,\pi) \subset L^2(-\pi,\pi) \rightarrow L^2(-\pi,\pi)$ of $J$ at $z_0$  is 
\begin{eqnarray}
dJ (z_0) z   = \lim_{\varepsilon \rightarrow 0} \frac{J\left(z_0+ \varepsilon z\right) - J\left(z_0\right)}{\varepsilon}=  \frac{\partial}{\partial x}\left(z_0 z\right).
\label{J'}
\end{eqnarray}
The linearized KS equation at $z_0$ is 
\begin{eqnarray}
\dot{z}=\left(A- dJ (z_0)  \right) z . 
\label{linear_ks_general}
\end{eqnarray}
The nonlinearity in the KS equation is not continuous and so the results for quasilinear systems do not apply. It also does not satisfy the assumptions of Theorem \ref{temam} as the linear operator in the KS equation is not negative and also  the nonlinear operator does not satisfy the assumptions. 
However, the $C_0$-semigroup $S\left(t\right)$ is Fr\'{e}chet differentiable at any $z_0 \in L^2(-\pi,\pi)$ and the derivative is the $C_0$-semigroup generated by the linearized KS equation at $z_0$. 
A similar result was shown in \cite{Temam-1988} but an additional assumption was required in the proof. 

\begin{thm}\cite{Rasha-2013,Rasha-2013c} Consider the nonlinear KS equation (\ref{cauchy-ks}). The nonlinear semigroup $S\left(t\right)$ is Fr\'{e}chet differentiable at every $z_0 \in L^2(-\pi,\pi)$. Moreover, indicating the  Fr\'{e}chet  derivative of $S$ by $T$,  $A+dJ( z_0) $ is the generator of $T$. 
\label{ks-frechet}
\end{thm}

The KS equation has an infinite number of equilibrium points. In particular, any constant function is an equilibrium to the KS equation.
%and there are an infinite number of other functions that are also equilibrium solutions but the explicit form is unknown \cite{Lan-2008, Michelson-1986}.% \achg{\st{for these equilibria although they can be approximated using Fourier truncated series. In this paper, we will} }
Define the closed invariant set of constant equilibria
\begin{eqnarray}
Z_e=\left\{ z_e: \ z_e \textup{ is a constant function} \right\} \subset  L^2(-\pi,\pi) . 
\label{equilibrium-set}
\end{eqnarray}

%Let $z_e$ be a constant function that does not depend on $x$, then $z_e$ is an equilibrium solution to the KS equation. The G\^{a}teaux derivative of the nonlinear operator $J$ at $z_e$ is 
%\begin{eqnarray*}
%dJ\left(z_e\right) z= z_e \frac{\partial z}{\partial x},
%\end{eqnarray*}
%and the linearized KS equation around a constant function $z_e$ is
%\begin{eqnarray}
%\dot{z}=Az - z_e  \frac{\partial z}{\partial x}. 
%\label{linear_ks}
%\end{eqnarray}

\begin{thm} 
Consider the KS equation (\ref{ks}). If the instability parameter $\nu>1$, then any $z_e \in Z_e$ is locally exponentially stable. If the instability parameter $\nu <1$, then the KS equation is unstable.
\end{thm}

\noindent 
\textbf{Proof.}
%For details, see \ref{Rasha-2013c}
Let $z_e$ be any constant equilibrium.
The operator $\left(A- dJ (z_e)   \right)$ is a Riesz-spectral operator with eigenvalues $\lambda_n = n^2 (1-\nu n^2 ) -i n z_e$, where $n\in \mathbb{Z}$ 
%and $z_e$ is a constant function that does not depend on $x$, and the corresponding eigenvectors are $\phi_n\left(x\right) = \frac{1}{\sqrt{2 \pi}} 
\cite[Theorem 5.2.1]{Rasha-2013}.  
Since $\left(A- dJ (z_e)  \right)$ is a Riesz-spectral operator, \cite[Theorem 2.3.5 c]{Curtain-1995} the spectrum determined growth assumption holds. The growth bound of the semigroup generated by $A+dJ (z_e) $ is determined by the supremum of the real part of the eigenvalues. 
%\begin{eqnarray*}
%\omega_0 = \sup_{n\in \mathbb{Z}} Re \{ \lambda_n \}.
%\end{eqnarray*}
Hence, if $\nu >1$, then all the eigenvalues of the linearized KS equation at a constant equilibrium have strictly negative real part, which results in a stable linearized system and if $\nu < 1$, then the linearized system is unstable. $\Box$

%The number of unstable eigenvalues  $N$ is the smallest integer such that 
%\begin{eqnarray}
%N > \sqrt{\frac{1}{\nu}}.
%\label{N}
%\end{eqnarray}

\section{CONCLUSIONS}

The application of Lyapunov's indirect method requires the $C_0$-semigroup of the nonlinear system to be Fr\'{e}chet differentiable. %\st{plays an important role in analyzing the stability using Lyapunov indirect method.} 
If the system linearized around an equilibrium is unstable or exponentially stable, then the equilibrium to the nonlinear system is unstable or  locally exponentially stable, respectively. If the  linearized system is only asymptotically stable, then Lyapunov's indirect method provides no conclusion about the stability of the equilibrium of the nonlinear system.%\achg{\sout{can be obtained}}
 This was illustrated by Example~\ref{example-Hans-Zwart}. %In addition, if the equilibrium \achg{\st{solution}} to the linearized system is unstable, then the nonlinear system is also unstable.}% A \achg{counterexample} by Hans Zwart was presented (Example \ref{example-Hans-Zwart}).

The key to using Lyapunov's Indirect Method is showing that the linearized generator corresponds to the generator of the  Fr\'{e}chet  derivative of the original semigroup. Since the generator is generally unbounded, this is not straightforward. Fairly complete results are available for quasi-linear systems, however.

Partial differential equations that are not quasi-linear are also considered. The only general result is Theorem \ref{temam}. There are many partial differential  equations that do not fit this class.  For example, the  Kuramoto-Sivashinsky equation  is not quasi-linear and does not satisfy the assumptions of Theorem \ref{temam}. However, it has a   Fr\'{e}chet  differentiable semigroup. Results using Lyapunov's Indirect Method  for the stability of this equation  are reviewed. 

Further research to establish general results on differentiability of semgroups and their relation to the linearized generator is needed in order to expand the applicability of Lyapunov's Indirect Method. 

\addtolength{\textheight}{-12cm}   % This command serves to balance the column lengths

%%%%%%%%%%%%%%%%%%%%%%%%%%%%%%%%%%%%%%%%%%%%%%%%%%%%%%%%%%%%%%%%%%%%%%%%%%%%%%%%

%%%%%%%%%%%%%%%%%%%%%%%%%%%%%%%%%%%%%%%%%%%%%%%%%%%%%%%%%%%%%%%%%%%%%%%%%%%%%%%%

%%%%%%%%%%%%%%%%%%%%%%%%%%%%%%%%%%%%%%%%%%%%%%%%%%%%%%%%%%%%%%%%%%%%%%%%%%%%%%%%
%\section*{APPENDIX}


%\section*{ACKNOWLEDGMENT}

%The preferred spelling of the word ÒacknowledgmentÓ in America is without an ÒeÓ after the ÒgÓ. Avoid the stilted expression, ÒOne of us (R. B. G.) thanks . . .Ó  Instead, try ÒR. B. G. thanksÓ. Put sponsor acknowledgments in the unnumbered footnote on the first page.

%%%%%%%%%%%%%%%%%%%%%%%%%%%%%%%%%%%%%%%%%%%%%%%%%%%%%%%%%%%%%%%%%%%%%%%%%%%%%%%%

%\rchg{I doubled check the bibstyle and I could not find another style for this conf. Could you check again?}
\bibliographystyle{IEEEtran}
\bibliography{references}

%\end{thebibliography}

\end{document}